\documentstyle{amsppt}
\magnification=1200
\NoBlackBoxes
\hsize=13cm
\def\C{\Bbb C}
\def\P{\Bbb P}

\def\ls{\vskip.25in}

\def\ms{\vskip.1in}
\def\Q{\Bbb Q}
\def\L{\Cal L}
\def\O{\Cal O}
\def\y{\bar{y}}
\def\bt{\boxtimes}
\null
\vskip 2cm
\centerline{\bf{SUBVARIETIES OF GENERIC HYPERSURFACES IN ANY VARIETY}}

\vskip 3cm
\baselineskip=12pt
\hskip .4cm {\bf LUCA CHIANTINI} \hskip 3cm {\bf ANGELO FELICE LOPEZ}
\vskip .3cm
\hskip .05cmDipartimento di Matematica \hskip 2.56cm
\ Dipartimento di Matematica

\hskip .7cm Universit\`a di Siena \hskip 3.9cm Universit\`a di Roma Tre

\hskip .65cm Via del Capitano 15 \hskip 3.2cm Largo San Leonardo Murialdo 1

\hskip .73cm 53100 Siena, Italy \hskip 4.5cm  00146 Roma, Italy 

\hskip .13cm e-mail {\tt chiantini\@unisi.it} \hskip 2.2cm e-mail 
{\tt lopez\@matrm3.mat.uniroma3.it}

\vskip .3cm
\centerline{AND}
\vskip .3cm
\centerline{\bf ZIV RAN}

\centerline{Department of Mathematics}

\centerline{University of California Riverside}

\centerline{Riverside, CA 92521, USA}

\centerline{e-mail {\tt{ziv\@math.ucr.edu}}}

\vfill\eject

{\bf {0. Introduction}}\ls
The geometry of a desingularization $Y^m$ of an arbitrary subvariety
of a generic hypersurface $X^n$ in an ambient variety $W$ (e.g. $W = \P^{n+1}$)
has received much attention over the past decade or so. Clemens [5] has
proved that for $m =1 , n = 2, W = \P^3$ and $X$ of degree $d$, $Y$ has genus
$g \geq 1 + d(d-5)/2$; Xu [13], [14] improved this to $g \geq d(d-3)/2 - 2$ for
$d \geq 5$ and showed that if equality holds for $d \geq 6$ then $Y$ is
planar; he also gave some lower bounds on the geometric genus
$p_g(Y)$ in case $m=n-1$. Voisin [10], [11] proved that, for $X$ of degree $d$ in
$\P^{n+1}, n\geq 3, m\leq n-2$ then
$p_g(Y)>0$ if $d\geq 2n+1-m$ and $K_Y$ separates generic points if
$d\geq 2n+2-m$ (see also [4], [1]).
For $X^n$ a generic complete intersection of type $(d_1,...,d_k)$ in any
smooth polarized $(n+k)$-fold $M$, Ein [6], [7] proved that $p_g(Y)>0$ if
$d_1+...+d_k\geq 2n+k-m+1$ and $Y$ is of general type if
$d_1+...+d_k\geq 2n+k-m+2.$
\par
In [2] the first two authors applied the classical method of
focal loci of families as in C. Segre [8] and Ciliberto-Sernesi [3]
(really normal bundle considerations) to give a new proof of one
of the main results of Xu in [13], [14] and to extend the lower bounds
for the genus in the cases of general surfaces in a component of
the Noether-Lefschetz locus in $\P^3$ and of general projectively
Cohen-Macaulay surfaces in $\P^4$.
\par
In this paper we introduce some notions of 'filling' subvarieties
and use them to prove two new results which {\it {inter alia}}
give another perspective on the above results, especially the genus
bounds. The general philosophy is as usual that as $Y$ moves with
$X$, sections of $K_Y$ (or some twist) can be produced through
differentiation; but here this is implemented by exploiting and
extending an elementary but perhaps surprising technique in the
spirit of classical projective geometry which goes back to [2]
and which gives a useful lower bound, depending
on the dimension of the projective span of $Y$,
on the number of independent
sections of $K_Y$ produced by the differentiation process.
\par
As to our results, in Theorem 1 below we extend and refine in the
above sense most of the above-quoted results (not including
Voisin's), by giving a lower bound on the number of sections of a
certain twist of $K_Y$ involving $K_W$ (for a recent extension, including
the proof of a conjecture of Clemens, see [15]).
Our second result (Theorem 2), based on the notion of
'$r$-filling' subvariety, indicates an apparently new and
unexpected direction as it deals with some higher-order tensors on
$Y$, manifested in the form of an effective divisor on a Cartesian
product $Y^r$ having certain vanishing order on a diagonal locus
as well as on a 'double point' locus associated to the map $Y\to
X$. As one application, we conclude a lower bound on the number of
quadrics (and higher-degree hypersurfaces) containing certain
'adjoint- type' projective images of $Y$ (and even on the
dimension of the kernel of certain 'symmetric Gaussian' maps)
(Corollary 2.1 below). Our feeling, however, is that the latter is
only the tip of an iceberg, and we hope to explore further in this
direction in the future.\ls
\vfill\eject

{\bf {1. Filling subvarieties}}\ls

Fix an $(n+1)$-fold $W$ (with isolated singularities)
and a free linear system $\L\subseteq H^0(L)$ for some line-bundle $L$, with
associated map $\varphi_{\L}$ to a projective space; set $\L_d$ for the
image of
$\L^{\otimes d}$ under the multiplication map
$H^0(L)^{\otimes d}\to H^0(L^d)$, and let
$X\in\L_d$ be a generic member. We consider generically finite maps
$$f:Y\to X\subset W$$
from a smooth $m$-fold. Such a map is said to be {\it{filling}} if
$(f,Y,X)$
deforms in a family $\{(f_t,Y_t,X_t):t\in T\}$ such that $X_t\in\L_d$ for all
$t$ ($X_t$ not fixed),
the natural map $T\to\L_d$ is generically finite, dominant and
$\bigcup_{t\in T}f_t(Y_t)$ is dense in $W$.
$f$ is said to be $(1,b)$-filling if, in addition, $(f,Y,X)$ moves
in a family $\{(f_s,Y_s,X):s\in S\}$ with $X$ fixed such that
$\bigcup_{s\in S}f_s(Y_s)$ is at least $(m+b)$-dimensional.
For instance,
if $W$ is homogeneous
then any $f$ is filling.\par

Our main result in the filling case is the following
which, though a special case of Theorem 2 below, we have
chosen to state and prove separately in order to make the
argument easier to follow.\ls
\proclaim{Theorem 1}
Let $X\in\L_d$ be generic, $f:Y\to X$ a
$(1,b)$-filling, generically finite map from a smooth irreducible $m$-fold,
and suppose $d \geq n-m-b$, $b\leq n-m$, $1\leq m\leq n-1$, and
that $\varphi_{\L}(f(Y))$ spans a $\P^{p+1}$.
Then
$$h^0(K_Y - (d-n+m+b)f^*L - f^*K_W)\geq 1+p(n-m-b). \tag 1 $$
If equality holds in (1) and $n-m-b\geq 1$,
then $f$ is not $(1,b+1)$-filling; in particular, if in addition $b=0$
then $f$ has no deformations which move $f(Y)$ with $X$ fixed.
\endproclaim

Before proving the Theorem, let us give a sampling of examples
and applications, beginning with the most popular case
$W=\P^{n+1}$. Set $d=2n+2-m+j$ and let $b=0$ unless otherwise mentioned.
Then we get:
$$h^0(K_Y(-j))\geq 1+p(n-m),\tag 2 $$
with equality only if $f$ has no deformations which move $f(Y)$ with $X$ fixed.\par
If $j\geq 0$ this bound implies a similar one for $K_Y$ itself, while Voisin
shows in this case that $\varphi_{K_Y}$ is generically 1-1 provided
$m\leq n-2$.\par
For $m=n-1$ and $j>0$, Xu gets a bound on $p_g$ while
the above shows $h^0(K_Y(-j))\geq 1+p$, which implies some
bounds on $p_g$ somewhat like Xu's; for a very lazy such bound,
note that our bound implies $p_g(Y)\geq 1+h^0(\O_Y(j))$
and the latter may be estimated using Koszul, assuming (as we may)
that $f(Y)$ is a $(d,e)$ complete intersection, yielding
$$p_g(Y)\geq 1+\binom{n+1+j}{j}-\binom{n+1+j-d}{j}-\binom{n+1+j-e}{j}
+\binom{n+1+j-d-e}{j}.\tag 3 $$
Perhaps better bounds can be obtained by estimating more
carefully the relevant multiplication map.\par
For $m=1$, i.e. $Y$ a curve,
and $b\geq 0$,
the basic estimate is
$$h^0(K_Y(-j))\geq 1+p(n-1-b),\tag 4 $$
with equality only if $f(Y)$ is not $(1,b+1)$-filling.
For instance,
a genus-3 nonplanar curve on a generic quintic surface is rigid
(we do not know if such curves exist); on
a generic sextic 3-fold $X$ any nonrigid curve $f(Y)$ must have
genus at least $1+p=$dim span$f(Y)$, and curves for which equality
holds don't fill up $X$
(note that 4-tangent plane sections have genus 6 and do fill up $X$);
the genus of a curve on a generic septic
3-fold is at least $1+2p$ (note that
a 6-tangent-plane section has genus 9).
\par
If $W$ is an abelian variety and $L$ is twice a principal
polarization (well-known to be free), and again taking $b=0$, note
that $p>0$ if $m>1$, so if in addition $d\geq n-m$, we conclude
that $h^0(K_Y)>1$; in particular \proclaim {} $X$ does not contain
any (translated) abelian subvariety of dimension
$m>1$.\endproclaim If $d>n-m, m\geq 1,$ we conclude as a special
case that \proclaim{} $X$ contains only subvarieties of general
type; in particular it does not contain any translated abelian
subvariety of dimension $m\geq 1$.\endproclaim
In Ein's general
situation, if $(M,L)$ is a smooth polarized $(n+k)$-fold and $L$
is assumed base-point free, we may take as $W$ any smooth complete
intersection of type $(d_2,...,d_k)$. By adjunction,
$K_W=(K_M+(d_2+...+d_k)L)|_W$. Letting $X$ be a generic member of
$H^0(L^d)|_W$, Theorem 1 yields \proclaim{}
$K_Y-(f^*K_M+(d_1+...+d_k-n+m+b)f^*L)$ is effective.\endproclaim

Now by Kodaira vanishing and Riemann-Roch, $K_M+(n+k+1)L$ is effective
(else the degree-$(n+k)$ polynomial $\chi (K_M+tL)$ would vanish at the
$n+k+1$ points $1,...,n+k+1$); moreover it is well-known [9]
that $K_M+(n+k)L$ is effective unless $M=\P^{n+k}$.
Note that as $f$ is filling, its image may be assumed not
contained in any fixed divisor, hence $f^*$ preserves effectiveness.
Thus Ein's result that $p_g(Y)>0$ if
$d_1+...+d_k\geq 2n+k-m+1$
follows, assuming $f$ is filling;
with Theorem 1 we conclude further, in this case, that
\proclaim{} If $(M,L)\not\simeq (\P^{n+k},\Cal O(1))$ then $p_g(Y)>0$
provided $d_1+...+d_k\geq 2n+k-m$.\endproclaim
If $f$ is not filling,
then by considering a desingularization of the subvariety of $M$ filled
up by $f(Y)$ we actually get a better bound. Note that Sommese [9]
has classified the polarized pairs $(M,L)$ such that $K_M+(n+k-2)L$ is not
$\Q$-effective, so excluding those and taking for simplicity $k=1$, for any
$f:Y\to X\subset M$ filling, we have that the Kodaira dimension of $Y$
is nonnegative (resp. $Y$ is of general type) provided $d\geq 2n-1-m$
(resp. $2n-m$).
Moreover note that in fact Ein proves $h^0(K_Y-(d+m-2n-2+b)f^*L) \not= 0$
when $f$ is $(1,b)$-filling, hence
in Theorem 1 we improve Ein's result in several directions: we give a
lower bound on the
dimension that depends upon how much $\varphi_{\L}(f(Y))$ spans, and we give
information on the extremal cases;
the twist is strictly better unless $M=\P^{n+k}$; and now $W$ can be any
variety with isolated singularities and not only a general complete
intersection in $M$.\smallskip

\demo {Proof of Theorem 1} Let $f:Y\to X\subset W$ be a generic member of
a $(1,b)$-filling family as above and let $N_f, N_{f/X}$ denote the
normal sheaves.
Note that in the definition of filling, there is no loss of generality
in assuming that the map $T\to\L_d$ is quasi-finite and unramified at the
point
corresponding to $X$; likewise in the definition of $(1,b)$-filling there
is no loss of generality in assuming that the natural map
$\coprod Y_s\to X$ has (differential) rank at least $m+b$ at a generic
point of $Y$.
 Denote by $N^0$ the subsheaf of $N_{f/X}$ generated
by $b$ generic global sections. Thus $N^0$ has rank
$b$ and $c_1(N^0)$
is effective. Now set $$N_{f/X}'=N_{f/X}/N^0, N_f'=N_f/N^0,$$
$$D=(\bigwedge^{n-m+1-b}N_f')^{**}.\tag 5$$ Thus
$$c_1(N_f)=K_Y-f^*K_W=D+D', D' {\text{effective}}.\tag 6$$
(In fact $D'=0$ iff $f$ is unramified in codimension 1.) Now
consider the exact
sequences
$$ 0\to N_{f/X}\to N_f\to f^*(L^d)\to 0, $$
$$ 0\to N_{f/X}'\to N_f'\to f^*(L^d)\to 0.$$

Our genericity assumption on $X$ implies that any infinitesimal
deformation of $X$ in $\L_d$ carries an infinitesimal deformation
of $f$, hence the natural map $\L_d\to H^0(f^*L^d)$ admits a
lifting $$\phi :\L_d\to H^0(N_f),\tag 7$$ whence sheaf maps $$\psi
:\L_d\otimes \Cal O_Y\to N_f,\tag 8$$ $$\psi' :\L_d\otimes \Cal
O_Y\to N_f'\tag 9$$ and pointwise maps at a general point $y\in Y$
$$\psi_y:\L_d\to N_{f,y}=\C^{n+1-m},\tag 10$$ $$\psi_y':\L_d\to
N_{f,y}'=\C^{n+1-m-b}.\tag 11$$ Our filling hypothesis implies
that $f(Y)$ moves, filling up $W$, hence that the map $\psi_y$ is
surjective, hence so is $\psi_y'$.\par The following technical
Lemma will easily imply Theorem 1. We will prove later a  more
general version, useful also for the applications in \S 2.

\proclaim{Lemma 1.1} Assumptions as in Theorem 1. Let $k$ be an
integer such that $0\le k\le d$ and let $$\Psi: \Cal
L_d\to\C^{k+1}$$ be a surjection of vector spaces. Then for a
general choice of elements $h_{k+1},\dots, h_{d}\in\Cal L$ and for
general subsets $Y_1,\dots,Y_k \subset Y$ each of cardinality
$p={\text{dim span}}\varphi_{\L}(f(Y))-1$,
the restriction of $\Psi$ to the subspace $\Cal L(-Y_1)\cdots \Cal
L(-Y_k)h_{k+1}\cdots h_d$ surjects.
\endproclaim

Let's first see that Lemma 1.1 implies Theorem 1. We apply it to
the map $\Psi=\psi_y'$ above with $k=n-m-b$. It yields the
existence of a generic subset $S\subset Y$ of cardinality
$p(n-m-b)$ such that $\psi_y'|_{\L_k(-S)h_{k+1}\cdots h_d}$ is
surjective. We conclude that the restriction of $\psi_y'$ on a
subbundle of the form $$\L_{(n-m-b)}(-S)(\prod_{a=1}^{d-(n-m-b)}
h_a)\otimes \O_{Y}$$ is generically surjective. Consequently, for
a generic $(n-m+1-b)$-dimensional subspace $V\subset\L_{(n-m-b)}$
we get an injective generically surjective map $\rho : V(\prod
h_a)\otimes\O_{Y}\to N_{f}'$ which evidently drops rank on $(\prod
f^*(h_a))_0$, as well as on $S$.\par Since
$\bigwedge^{n-m+1-b}\rho$ yields a section of
$\bigwedge^{n-m+1-b}N_{f}^{'**}$ vanishing on a generically placed
divisor of class $(d-(n-m-b))f^*L$ as well as on $p(n-m-b)$ many
generic points, it follows that
$$h^0(c_1(N_{f}^{'**})-(d-(n-m-b))f^*L)\geq 1+p(n-m-b)$$ and
Theorem 1 follows in light of Eq. (5)(6).\qed\enddemo

Now Lemma 1.1 is the special case $r=1$ of the following result
that we are going to use in its full generality in the next
section.

\proclaim{Lemma 1.2} Assumptions as in Theorem 1. Fix integers
$r,k$ with $0\le r-1\le k\le d$. Let $\Cal L_d\to \C^N$ be a
linear map of vector spaces such that for $y_1,\dots,y_{r-1}$
general points of $Y$, the restriction to $\Cal
L_d(-y_1-\dots-y_{r-1})$ induces a surjection
$$\Psi:\Cal
L_d(-y_1-\cdots-y_{r-1})\to\C^{k+1}.$$
Then for a general choice
of elements $h_{k+1},\dots, h_d\in\Cal L$ and for general subsets
$Y_1,\dots,Y_k \subset Y$ each of cardinality $p$, with $Y_i\ni
y_i$, $i=1,\dots,r-1$,  the restriction of $\Psi$ to the subspace
$\Cal L(-Y_1)\cdots \Cal L(-Y_k)h_{k+1}\cdots h_d$ surjects.
\endproclaim
\demo{Proof} Fix an integer $j$ such that $0 \le j \le k$. We say
that $U\subset \Cal L_d(-y_1-\cdots -y_{r-1})$ is a subspace {\it
of type $j$} if $$U=U(Y_1,\dots,Y_j,h_{j+1},\dots,h_d):= \Cal
L(-Y_1)\cdots\Cal L(-Y_j)h_{j+1}\cdots h_d$$ for some subsets
$Y_1,\dots,Y_j \subset Y$ of cardinality $p$ and elements
$h_{j+1},\dots,h_d$ in $\Cal L$, such that there exists an integer
$q$ with $0 \le q \le \min\{j,r-1\}$ satisfying $y_1\in Y_1,\dots,
y_q\in Y_q, h_{j+1}\in \Cal L(-y_{q+1}),\dots, h_{r+j-1-q}\in \Cal
L(-y_{r-1})$. Given $U$, the maximal such integer $q$ will be
called the {\it index} of $U$. Observe that a general subspace of
type $j<d$ and index $q$ is contained in a general subspace of
type $j+1$, index $\ge q$. Here of course the word 'general' means: in
a general choice of the linear forms $h_i$ and the subsets $Y_i$,
subject only to the restrictions determined by the points
$y_1,\dots, y_{r-1}$ and the index $q$.
Let us now record the following
\proclaim {Claim 1} If for a general
subspace of type $j$ and index $q < \min\{j,r-1\}$ one has
$\dim\Psi(U)=x$, then for a general subspace $U'$ of type $j$ and
index $q+1$ one has $\dim\Psi(U')\ge x$.\endproclaim
\demo {Proof of Claim 1} Observe indeed that the points $y_1,\dots,
y_{r-1}$ are chosen generically. On the other hand, if $U=\Cal
L(-Y_1)\cdots\Cal L(-Y_j)h_{j+1}\cdots h_d$ is of type $j$ and
index $q$ for $y_1,\dots,y_{r-1}$, then since $Y_{q+1}$ is
general, without any loss of generality we may replace $y_{q+1}$
with some point of $Y_{q+1}$ and similarly replace the $h$'s;
hence $U$ becomes of index $q+1$. This proves Claim 1.\enddemo
The proof of Lemma 1.2 is based on an inductive construction of
general subspaces $U_0, U_1,\dots, U_k$ of $\Cal
L_d(-y_1-\cdots-y_{r-1})$, such that for all $j$, $U_j$ is of type
$j$ and the image under $\Psi$ has dimension $\ge
j+1$.\par\noindent Since $\Cal L_d(-y_1-\cdots-y_{r-1})$ is
generated by monomials $h_1\cdots h_d$ and $\Psi$ surjects, then
for a general choice of $U= \C h_1 \cdots h_d$ of type $0$, we get
$\dim\Psi(U)=1$. Assume now that, for some $j<k$, we have
constructed a general subspace $U$ of type $j$ and maximal index
$q=\min\{j,r-1\}$ whose image under $\Psi$ is of dimension $\ge
j+1$. If this dimension is actually bigger than $j+1$, then a
general subspace of type $j+1$ and maximal index has image of
dimension $\ge j+2$ and the induction goes on. Hence we can assume
that for a general $U$ of type $j$ we have $\dim\Psi(U)=j+1$.
Since $j<k$, and $\Psi$ is surjective we know that for a general
monomial $M\notin U$ we have $\Psi(M)\notin\Psi(U)$.\par
\proclaim
{Claim 2} For any $j=0,\dots,k$, let $U=\Cal
L(-Y_1)\cdots\Cal L(-Y_j)h_{j+1}\cdots h_d$ be a general subspace
of type $j$, index $q=\min\{j,r-1\}$ and let $M'=h'_1\cdots h'_d$
be a general monomial in $\Cal L_d(-y_1-\cdots -y_{r-1})$. Then
there exists a chain of subspaces $U^0=U,U^1,\dots,U^a$ of type
$j$, such that $M'\in U^a$ and each consecutive pair $U^i,U^{i+1}$
is contained in the same subspace of type $j+1$ and index $\ge q$.
Furthermore as $U$,$M'$ move generically, we may assume that any
subspace of the chain is general.\endproclaim
\demo {Proof of Claim 2} Without
any loss of generality, we may assume that $h'_i\in\Cal L(-y_i)$,
for all $i=1,\dots,r-1$. Let us now construct a chain as above,
which links $U$ and the subspace $\Cal L(-Y_1)\cdots\Cal
L(-Y_j)h'_{j+1}h_{j+2}\cdots h_d$. First suppose $j\ge r-1$. Pick
a general subset $\{A_1,\dots, A_p\} \subset Y$, over which
$f^*h_{j+1}$ vanishes, and a general subset $\{B_1,\dots,B_p \}
\subset Y$ contained in the locus where $f^*h'_{j+1}=0$.  Take,
for $c=1,\dots, p$, a general element $z_c\in\Cal
L(-B_1-\dots-B_c-A_c-\dots -A_p)$, which exists by our definition
of $p$. Then put $z_0=h_1$, $z_{p+1}=h'_1$ and define for all
$c=0,\dots,p+1$, $U^c=\Cal L(-Y_1)\cdots\Cal
L(-Y_j)z_ch_{j+2}\cdots h_d$. It is clear that $U^c,U^{c+1}$ both
are contained in $\Cal L(-Y_1)\cdots\Cal L(-Y_j)\Cal
L(-B_1-\dots-B_c-A_{c+1}-\dots -A_p)h_{j+2}\cdots h_d$, which is
in fact  general of type $j+1$, for $h_{j+1},h'_{j+1}$ and the
points $A_i$'s, $B_i$'s are general. When $j<r-1$, i.e. when both
$h_{j+1},h'_{j+1}$ lie in $\Cal L(-y_{j+1})$, we argue as above,
except that we take $B_p=A_p=y_{j+1}$. We repeat then the previous
procedure, considering $h_{j+2},h'_{j+2}$ instead of $h_{j+1},
h'_{j+1}$ and so on. We get finally a chain with the required
properties, that links $U$ with the subspace $\Cal
L(-Y_1)\cdots\Cal L(-Y_j)h'_{j+1}\cdots h'_d$. Therefore we may
replace, from now on,  $U$ with the subspace $\Cal
L(-Y_1)\cdots\Cal L(-Y_j)h'_{j+1}\cdots h'_d$. Choose now a
general subset $Y'\subset Y$ of cardinality $p$, containing
$y_{j+1}$ if $j+1<r$ and such that $h'_{j+1}\in\Cal L(-Y')$.
Choose also a general $h\in \Cal L(-Y_j)$. Define $U^*=\Cal
L(-Y_1)\cdots\Cal L(-Y_{j-1})\Cal L(-Y')hh'_{j+2}\cdots h'_d$. By
our generality assumptions on $Y'$ and $h$, also $U^*$ is a
general subspace of type $j$ and maximal index and both $U$ and
$U^*$ are contained in $\Cal L(-Y_1)\cdots\Cal L(-Y_j)\Cal
L(-Y')h'_{j+2}\cdots h'_d$, which is of type $j+1$. Now use the
procedure at the beginning of the proof, acting on $h$ and $h'_j$,
to find a chain of the required type which links $U^*$  with the
subspace $\Cal L(-Y_1)\cdots\Cal L(-Y_{j-1})\Cal
L(-Y')h'_jh'_{j+2}\cdots h'_d$. At the end of this procedure, we
see that we may link $U$ with the new subspace $\Cal
L(-Y_1)\cdots\Cal L(-Y_{j-1})\Cal L(-Y')h'_jh'_{j+2}\cdots h'_d$;
notice that renumbering $Y'$ as $Y_j$ and swapping $h'_j$ with
$h_{j+1}$, we linked $U$ with the subspace $\Cal L(-Y_1)\cdots\Cal
L(-Y_j)h'_{j+1}\cdots h'_d$, but with the good news that now
$h'_j\in\Cal L(-Y_j)$. Repeat this last construction, using $h\in
L(-Y_{j-1})$ and so on. At the very end we see that we constructed
a chain with the required properties, linking our original $U$
with a subspace of type $\Cal L(-Y_1)\cdots\Cal
L(-Y_j)h'_{j+1}\cdots h'_d$ for which $h'_1\in\Cal L(-Y_1),\dots,
h'_j\in\Cal L(-Y_j)$. This is our $U^a$ and Claim 2 is
proved.\enddemo
 Let us go back to the proof of Lemma 1.2. By Claim
2, we get a chain as above of subspaces of type $j$ joining $U$
with $M$; all the elements of the chain are general, so replacing
$U$ with some element of the chain, we may assume that $U,M$ both
belong to the same subspace $U'$ of type $j+1$. Since $U$ is
general, $U'$ is also general (of its index) and on the other hand
$\dim\Psi(U')>\dim\Psi(U)$. Therefore we have the inductive step.
\qed\enddemo

{\bf Remark.} The case $k=1$ is in [2]; the general case is not
much more difficult. Another proof of the above lemma can be found in 
[15]. \ls 

{\bf{2. $r$-filling subvarieties}}\ls We
continue with the notations of the previous section. The map $f$ is
said to be $r$-filling if $(f,Y,X)$ deforms in a family
$\{(f_t,Y_t,X_t):t\in T\}$ such that $X_t\in\L_d$ for all $t$ ($X_t$
not fixed), the natural map $T\to\L_d$ is generically finite, dominant and
$\bigcup_{t\in T}f_t(Y_t)^r$ is dense in the Cartesian product
$W^r$. $f$ is said to be $(r,b)$-filling if it is $r$-filling and
$(1,b)$-filling.
For instance,
\par
- if $W=\P^{n+1}$ and $f(Y)$ spans (at least) a $\P^{r-1}$, then
clearly $f$ is $r$-filling (since linearly independent $r$-tuples
are projectively equivalent);\par - if $W$ is a Grassmannian of
lines in $\P^a$ then, because all pairs of disjoint lines are
projectively equivalent, it follows that any $f$ is 2-filling
unless $f(Y)$ is contained in a (Pl\"ucker-) linear subvariety
$U\subset W$, and it is easy to see that any such $U$ must be of
the form either \{ lines through a point \} or \{ lines in a $\P^2
\}$.\smallskip The following result extends Theorem 1 to the
$r$-filling case:\ls \proclaim{Theorem 2} Let $X\in\L_d$ be
generic, $f:Y\to X$ an $(r,b)$-filling, generically finite map
from a smooth irreducible $m$-fold, and suppose $d \geq r(n-m-b)$,
$r\leq n-m+1-b$, $1\leq m\leq n-1$, and that $\varphi_{\L}(f(Y))$
spans a $\P^{p+1}$. Let $\Delta_Y^r\subset Y^r$ be the big
diagonal (= locus of nondistinct $r-$tuples) and $D_f^r\subset
Y^r$ the double locus of $f$ (=  closure of locus of distinct
$r-$tuples $(y_1,...,y_r)$ such that $f(y_1),...,f(y_r)$ are not
distinct). Then the linear system of $S_r$-invariant sections in
$H^0(K_Y-(d-r(n-m-b))f^*L-f^*K_W)^{\bt r})$ having multiplicity at
least $n-m+1$ on $\Delta_Y^r$ and multiplicity at least $n$ on
$D_f^r$ is of dimension at least $2-r+pr(n-m-b)$.\endproclaim
To
get an idea of the significance of this result for $r>1$ consider
now the case $r=2, n-m-b\geq 1$. Set
$M=K_Y-f^*K_W-(d-2(n-m-b))f^*L$, and assume $f$ is 2-filling
(which automatically holds if $W=\P^{n+1}$ in which case 
$M=K_Y(3n-2m-d-2b+2)$). The Theorem 2 yields
a space of symmetric elements $\tau\in H^0(M)^{\otimes 2}$
vanishing to order $n-m+1\geq 2$ on the diagonal. These in turn
yield, in particular, elements of $$I_2(M)=\ker ({\text {sym}}
^2(H^0(M))\to H^0(M^{\otimes 2}))$$ i.e. quadrics containing
$\varphi_M(Y)$ (in fact, if $n-m>1$ they belong to the kernel of a
certain 'higher-order symmetric Gaussian', cf. [12]); thus

\proclaim{Corollary 2.1} In the above situation, we have $$\dim
I_2(M)\geq 2p(n-m-b).$$\endproclaim Clearly similar results can be
stated for $r>2$; we shall leave this to the reader.\par
\demo{Proof of Theorem 2} We continue with the notations used in
the proof of Theorem 1. Now set $\y = \{y_1,..., y_{r-1}\}$ for
generic $y_1,..., y_{r-1}\in Y$. As before we have sheaf maps $\psi
:\L_d\otimes \Cal O_Y\to N_f,$ $\psi' :\L_d\otimes \Cal O_Y\to
N_f'$ and pointwise maps $\psi_y:\L_d\to N_{f,y}=\C^{n+1-m},
\psi_y':\L_d\to N_{f,y}'=\C^{n+1-m-b}, y\in Y$ general (cf.
(7)-(11)). Our $r$-filling hypothesis implies that $f(Y)$ moves,
filling up $W$ while fixing $f(y_1),...,f(y_{r-1})$.
 We claim that
this (and the hypothesis on $d$), implies that the restricted map
$$\Psi_y:\L_d(-\y )\to N_{f,y}\oplus N_{f/X,\y}$$ is surjective,
where $\L_d(-\y )$ denotes the set of elements of $\L_d$ vanishing
on $f(y_1),...,$\par \noindent $f(y_{r-1})$ and $N_{f/X,\y}
=\bigoplus\limits_{i=1}^{r-1}N_{f/X,y_i}$. 
Indeed this follows immediately by the 'snake lemma' from the fact
that we have a surjection $\L_d\to N_{f,y}\oplus N_{f,\y}$ ($r$-
fillingness) which induces an isomorphism
$L^{\otimes d}\otimes\O_{\y}\to N_{f,\y}/N_{f/X,\y}.$
Or more concretely, given any
vector of type $(v,0) \in N_{f,y}\oplus N_{f/X,\y}$ we can find an
element of $\L_d(-\y )$ inducing the deformation $v$ on $y$, hence
the vectors $(v,0)$ lie in the image of $\Psi_y$. Similarly given
any vector of type $(0, w_{i}) \in N_{f,y}\oplus N_{f/X,\y}$ with
$w_{i} \not= 0$ only on its i-th component, again by
r-fillingness, we can find a vector $w \in N_{f,y}$ such that
$(w,w_{i})$ lies in the image of $\Psi_y$. Therefore $\Psi_y$ is
surjective and likewise the map $\Psi_y':\L_d(-\y )\to
N_{f,y}'\oplus N'_{f/X,\y}$ is surjective  as well. Now consider
the Cartesian product $$f^r:Y^r\to X^r\subset W^r .$$ We have an
exact $S_r$-invariant sequence $$0\to N_{f^r/X^r}\to N_{f^r}\to
f^*(L^d)^{\boxplus r}\to 0$$ and $N_{f^r}= N_f^{\boxplus r}$ etc.
Likewise $$0\to N_{f^r/X^r}'\to N_{f^r}'\to f^*(L^d)^{\boxplus
r}\to 0, $$ $N'_{f^r}=(N'_f)^{\boxplus r}$. Now there is a natural 
$S_r-$invariant (i.e. symmetric) map
$$\L_d\to H^0((f^*L^d)^{\boxplus r})=\oplus H^0(f^*(L^d))$$ and by
the $r$-filling assumption this lifts to 
a symmetric map $\L_d\to H^0(N_{f^r})$, whence a map
$\phi^r:\L_d\to H^0(N_{f^r}')$ whose corresponding sheaf map
$\psi^r:\L_d\otimes \O_{Y^r}\to N_{f^r}'$ is generically
surjective and symmetric.\par
We now apply Lemma 1.2 to the map $\Psi=\Psi_y'$
above with $k=r(n-m-b)$. It yields the existence of a generic
(independently of $\y$) subset $S\subset Y$ of cardinality
$pr(n-m-b)-r+1$ such that $\Psi_y'|_{\L_k(-S)h_{k+1}\cdots h_d}$
is surjective. We conclude that the restriction of $\psi^r$ on a
subbundle of the form $$\L_{r(n-m-b)}(-S)(\prod_{a=1}^{d-r(n-m-b)}
h_a)\otimes \O_{Y^r}$$ is generically surjective. Consequently,
for a generic $r(n-m+1-b)$-dimensional subspace
$V\subset\L_{r(n-m-b)}$ we get an injective generically surjective map 
$$\rho^r : V(\prod h_a)\otimes\O_{Y^r}\to N_{f^r}'$$ 
and since the $h_i$ are chosen generally, independently
of the $y_j, \rho^r$ is symmetric as well. 
On the other hand $\rho^r$ evidently drops rank on $(\prod f^*(h_a))_0$, 
as well as on
slices $Y^{r-1}\times S$ (and their $S_r-$orbits).\par 
Now let $$b_r:Y_r\to Y^r$$ be the
blowing-up of the components of the big diagonal in some order,
with exceptional divisor $E_{Y,r}$, and likewise for $X,W$ (same
order). We have a natural map $f_r:Y_r\to X_r\subset W_r$, whence
an exact sequence $$0\to N_{f_r/X_r}\to N_{f_r}\to
f^*(L^d)^{\boxplus r}\to 0 .$$ Let
$N'_{f_r}=N_{f_r}/b_r^*(N^{0\boxplus r})$ and note that
$c_1(N_{f_r})=c_1(N_{f_r}')+A, A$ effective. Since $f_r$ moves
with $f$, there is as above a generically surjective map
$\psi_r:\L_d\otimes\O_{Y_r}\to N_{f_r}'$ whence a generic
isomorphism $\rho_r:V\otimes\O_{Y_r}\to N_{f_r}'$ compatible with
$\rho^r$. Now we compute $$c_1(N_{f_r})=b_r^*(K_Y-f^*K_W)^{\bt
r}+(m-1)E_{Y,r}-nf_r^*E_{W,r}$$ $$=b_r^*(K_Y-f^*K_W)^{\bt
r}-(n-m+1)E_{Y,r}-nR,$$ where $R$ is effective and contains the
blowup of some scheme supported on the locus of distinct
$r$-tuples mapped by $f$ to nondistinct ones. Since
$\bigwedge^{r(n-m+1-b)}\rho^r$ yields a symmetric section of
$\bigwedge^{r(n-m+1-b)}N_{f_r}^{'**}$ vanishing on a generically
placed symmetric divisor of class $(d-r(n-m-b))L^{\bt r}$ as well
as on $pr(n-m-b)-r+1$ many slices $Y^{r-1}\times$pt., pt.$\in S$,
the Theorem 2 follows by projecting back to $Y^r$.\qed\par 

{\bf Remark.} Rather than blow up the big diagonal we could equally
well blow up other 'diagonal strata', thus leading to analogues of
the Theorem 2 which we will let the reader state; though we do not
pursue them here, they might prove useful.\par
\enddemo
\vfill\eject

\centerline{\bf References}
\vskip .5cm

\item{[1]}  M.\ C.\ Chang and Z.\ Ran.\ Divisors on Generic hypersurfaces.\
{\it J.\ Differential Geom.\ \bf 38} (1993), 671-678.
\ms
\item{[2]} L.\ Chiantini and A.\ F.\ Lopez.\ Focal loci of families
and the genus of curves on surfaces.\ {\it Proc. Amer.\ Math.\ Soc.\ 
\bf 127} (1999), 3451-3459. 
\ms
\item{[3]} C.\ Ciliberto and E.\ Sernesi.\ Singularities of the theta
divisor and congruences of planes.\ {\it J.\ Alg.\ Geom.\ \bf 1}
(1992), 231-250.
\ms
\item{[4]} H.\ Clemens.\ Curves in generic hypersurfaces.\ {\it Ann.\ 
Sci.\ \'Ec.\ Norm.\ Sup.\ \bf 19} (1986), 629-636.
\ms
\item{[5]} H.\ Clemens, J.\ Koll\'ar and S.\ Mori.\ Higher-dimensional 
complex geometry.\ {\it Ast\'eris\-que 166} (1988).
\ms
\item{[6]} L.\ Ein.\ Subvarieties of generic complete intersections.\
{\it Invent.\ Math.\ \bf 94} (1988), 163-169.
\ms
\item{[7]} L.\ Ein.\ Subvarieties of generic complete intersections II.\
{\it Math.\ Ann.\ \bf 289} (1991), 465-471.
\ms
\item{[8]} C.\ Segre.\ Sui fuochi di $2^o$ ordine dei sistemi
infiniti di piani, e sulle curve iperspaziali con una doppia
infinit\`a di piani plurisecanti.\ {\it Atti R.\ Accad.\ Lincei
\bf 30} vol. 5, (1921), 67-71.
\ms
\item{[9]} A.\ Sommese.\ On the adjunction-theoretic structure of projective
varieties.\ In {\it Complex analysis and algebraic geometry (G\"ottingen, 1985)},
Lecture Notes in Math.\ vol.\ {\bf 1194} (Springer, Berlin-New York, 1986), 
175-213.
\ms
\item{[10]} C.\ Voisin.\ On a conjecture of Clemens on rational
curves on hypersurfaces.\ {\it J.\ Differential Geom.\ \bf 44} (1996), 200-213.
\ms
\item{[11]} C.\ Voisin.\ A correction: "On a conjecture of Clemens on rational 
curves on hypersurfaces".\ {\it J.\ Differential Geom.\ \bf 49}, 1998, 601-611.
\ms
\item{[12]} J.\ Wahl.\ Introduction to Gaussian maps on an algebraic curve.\
In {\it Complex projective geometry (Trieste, 1989/Bergen, 1989)}, London Math.\ 
Soc.\ Lecture Note Ser.\, {\bf 179} (Cambridge Univ.\ Press, 1992), 304-323.
\ms
\item{[13]} G.\ Xu.\ Subvarieties of general hypersurfaces in projective
space.\ {\it J.\ Differential Geom.\ \bf 39} (1994), 139-172.
\ms
\item{[14]} G.\ Xu.\ Divisors on generic complete intersections in projective 
space.\ {\it Trans.\ Amer.\ Math. Soc.\ \bf 348} (1996), 2725-2736.
\ms
\item{[15]} Z.\ Ran.\ On a conjecture of Clemens.\ {\it Preprint} math-ag 9911161
(1999).

\enddocument